\documentclass[11pt,reqno,a4paper]{article}

\usepackage[blocks]{authblk}
\usepackage{amsmath,amsthm,amstext,amscd,amssymb,euscript,url,mathtools,mathrsfs,euscript,dsfont,comment,epsfig,todonotes,relsize,scalerel}
\usepackage{hyperref}
\usepackage[shortlabels]{enumitem}

\definecolor{unbleu}{rgb}{0.03, 0.15, 0.4}

\hypersetup{
pdfborder = {0 0 0},
colorlinks,
linkcolor=unbleu,
citecolor=unbleu,
urlcolor=unbleu
}

\usepackage{xcolor}

\newcommand{\sx}{\scaleto{(x)}{6pt}}
\newcommand{\sy}{\scaleto{(y)}{6pt}}
\newcommand{\Z}{\mathds Z}
\newcommand{\R}{\mathds R}
\newcommand{\N}{\mathds N}

\newcommand{\E}{\mathds E}

\newcommand{\Zd}{\mathds Z^d}

\renewcommand{\phi}{\varphi}
\newcommand{\epsi}{\ensuremath{\epsilon}}
\newcommand{\La}{\ensuremath{\Lambda}}
\newcommand{\la}{\ensuremath{\Lambda}}
\newcommand{\si}{\ensuremath{\sigma}}

\newcommand{\pee}{\ensuremath{\mathbb{P}}}

\def\1{{\mathds 1}}

\newcommand{\dd}{\mathop{}\!\mathrm{d}}

\newtheorem{theorem}{{\small T}{\scriptsize HEOREM}}[section]
\newtheorem{corollary}{{\bf{\small C}{\scriptsize OROLLARY}}}[section]
\newtheorem{proposition}{{\bf{\small P}{\scriptsize ROPOSITION}}}[section]
\newtheorem{lemma}{{\bf{\small L}{\scriptsize EMMA}}}[section]
\newtheorem{remark}{{\bf{\small R}{\scriptsize EMARK}}}[section]
\newtheorem{definition}{{\bf{\small D}{\scriptsize EFINITION}}}[section]
\newtheorem{induction}{{\bf{\small I}{\scriptsize NDUCTIVE HYPOTHESIS}}}[section]

\renewenvironment{proof}[1]
{\noindent{{\bf{\small{ P}{\scriptsize ROOF}}}.}\hspace{0.1cm} #1} {$\;\qed$\newline}

\newcommand{\beq}{\begin{eqnarray}}
\newcommand{\eeq}{\end{eqnarray}}

\newcommand{\ba}{\begin{align*}}
\newcommand{\ea}{\end{align*}}

\newcommand{\be}{\begin{equation}}
\newcommand{\ee}{\end{equation}}

\newcommand{\bl}{\begin{lemma}}
\newcommand{\el}{\end{lemma}}

\newcommand{\br}{\begin{remark}}
\newcommand{\er}{\end{remark}}

\newcommand{\bt}{\begin{theorem}}
\newcommand{\et}{\end{theorem}}

\newcommand{\bd}{\begin{definition}}
\newcommand{\ed}{\end{definition}}

\newcommand{\bind}{\begin{induction}}
\newcommand{\eind}{\end{induction}}

\newcommand{\bp}{\begin{proposition}}
\newcommand{\ep}{\end{proposition}}

\newcommand{\bc}{\begin{corollary}}
\newcommand{\ec}{\end{corollary}}

\newcommand{\bpr}{\begin{proof}}
\newcommand{\epr}{\end{proof}}

\newcommand{\bi}{\begin{itemize}}
\newcommand{\ei}{\end{itemize}}

\newcommand{\ben}{\begin{enumerate}}
\newcommand{\een}{\end{enumerate}}

\newcommand{\caC}{{\mathscr C}}

\newcommand{\caF}{{\mathcal F}}

\newcommand{\caI}{{\mathcal I}}

\DeclareMathOperator{\ent}{\scaleto{s}{5.5pt}}

\newcommand{\sipr}{\{\si(t), t\geq 0\}}
\newcommand{\gcb}[1]{\mathrm{GCB}\big(#1\big)}

\DeclareMathOperator{\e}{\mathrm{e}}

\begin{document}

\title{Gaussian concentration bounds for \\ probabilistic cellular automata}

\author[1]{Jean-Ren\'e Chazottes
\thanks{Email: \texttt{jeanrene@cpht.polytechnique.fr}}
}
\author[2]{Frank Redig
\thanks{Email: \texttt{F.H.J.Redig@tudelft.nl}}
}

\author[2]{Edgardo Ugalde
\thanks{Email: \texttt{ugalde@ifisica.uaslp.mx}}
}

\affil[1]{{\small Centre de Physique Th\'eorique, CNRS, Ecole polytechnique, Institut Polytechnique de Paris, France}}
\affil[2]{{\small Institute of Applied Mathematics, Delft University of Technology, The Netherlands}}
\affil[3]{{\small Instituto de F\'{\i}sica,Universidad Aut\'onoma de San Luis Potos\'{\i}, M\'exico}}

\date{Dated: \today}

\maketitle

\begin{abstract}
We study lattice spin systems and analyze the evolution of Gaussian concentration bounds (GCB) under the action of probabilistic cellular automata (PCA), which serve as discrete-time analogues of Markovian spin-flip dynamics.
We establish the conservation of GCB and, in the high-noise regime, demonstrate that GCB holds for the unique stationary measure. Additionally, we prove the equivalence of GCB for the space-time measure and its spatial marginals in the case of contractive probabilistic cellular automata.
Furthermore, we explore the relationship between (non)-uniqueness and GCB in the context of space-time Gibbs measures for PCA and illustrate these results with examples.
\end{abstract}

\newpage


\section{Introduction}
In this paper we are interested in proving  concentration inequalities for certain measures far from product measures.
Concentration inequalities imply that if a function of many weakly dependent random variables only weakly depends on each of them individually, then the function concentrates around its expected value.
A key aspect of this concentration phenomenon is its non-asymptotic nature, distinguishing it from classical limit theorems, where results emerge only as the number of random variables tends to infinity. Recall that the three main types of limit theorems are the law of large numbers, the central limit theorem, and large deviations.
Another crucial advantage of concentration inequalities is their flexibility: they apply to arbitrarily defined functions of random variables, provided these functions are sufficiently regular (e.g., Lipschitz or bounded differences), whereas classical limit theorems primarily address sums of random variables.

Concentration inequalities have led to a paradigm shift in probability and statistics, extending their impact to discrete mathematics, geometry, and functional analysis (see, e.g., \cite{blm,dp,ledoux,vershynin,Wbook}).

For functions of independent random variables satisfying the bounded-differences property, a fundamental result is McDiarmid's inequality, which exhibits a Gaussian-type of concentration.
In the setting of Gibbs measures in lattice spin systems, McDiarmid's inequality can be seen as the high-temperature limit of a more general Gaussian concentration bound (GCB), which applies to all local functions. As the temperature tends to \(+\infty\), interactions become negligible, and the system approaches a product measure, recovering McDiarmid's result as a special case.  GCB indeed holds under high temperature conditions, such as the Dobrusshin uniqueness  \cite{kuelske} and its extensions, including complete analyticity \cite{dobshlos}. In finite-range lattice spin systems, complete analyticity is known to be equivalent to the log-Sobolev inequality \cite{stroock}, which, in turn, implies GCB \cite{ledoux}. Various applications are given in \cite{ccr}.

In this paper, we investigate the  GCB property  in the setting of  probabilistic cellular automata (PCA) -- both for their  stationary measures  and their
space-time Gibbs measures \cite{lms}.
PCA are a class of models where an infinite collection of cells, each positioned at the vertex of a lattice (here it will be $\Z^d$), evolve over discrete time steps according to stochastic dynamics. Each cell can occupy one of a
finite number of states. The evolution of these cells is governed by Markov chains, where the state of a cell at any given time depends probabilistically on the states of its neighboring cells at the previous
time step. This framework combines the spatial structure of cellular automata with the probabilistic transitions of Markov chains, enabling the study of complex systems with inherent randomness \cite{mm}.
More precisely, we study and partially answer the following key questions:

\begin{enumerate}
\item
Time evolution of GCB:  If a PCA starts from a measure  $\mu $ satisfying GCB, is this property preserved over time? How does the associated concentration constant evolve?
Does GCB hold in the limit for stationary measures? Furthermore, can GCB provide insights into the speed of convergence to equilibrium?
\item GCB and space-time measures:  How is the GCB property of a stationary measure related to the GCB property of the corresponding  space-time (Gibbs) measure ?
\item GCB and uniqueness:  How does the GCB property of a stationary measure relate to  uniqueness  and GCB of the corresponding stationary space-time Gibbs measure?
\end{enumerate}

These questions naturally extend those studied in the context of  spin-flip dynamics  in continuous time in \cite{ccr2}. They also relate to the connection between GCB and
uniqueness of  translation-invariant Gibbs measures  in lattice spin systems \cite{moles,cr}. Moreover, they fit into the broader framework of  space-time Gibbs measures
developed in \cite{lms}, where stationary (translation-invariant) measures are seen as restrictions of space-time Gibbs measures to temporal layers. The interplay between the
GCB property of such restrictions and that of the original measure is, in general, subtle. For instance, in the low-temperature Ising model, non-uniqueness and GCB of
restrictions can coexist, as we will illustrate with an example.

Beyond its theoretical significance, another motivation for studying GCB in this context is to establish concentration properties for  measures that are not available in an explicit form (e.g. via their finite dimensional distributions, or via a Gibbsian potential)  but which
can be obtained as a time evolution or as limiting stationary distributions. In  non-equilibrium systems where detailed balance or reversibility does not hold, explicit
descriptions of stationary measures are rare (e.g., in the form of a Gibbs measure). Consequently, their concentration properties cannot necessarily be inferred from standard high-
temperature criteria, making alternative approaches necessary.

{\em Organization of the Paper}.
In Section \ref{pcasect}, we introduce the fundamental definitions and concepts.
In Section \ref{pcaev}, we prove the conservation of Gaussian concentration under PCA evolution and establish the equivalence between space-time GCB and concentration
of spatial marginals for  contractive PCA. Additionally, we provide a result on the speed of relaxation to equilibrium under a weaker contraction property (less restrictive than
exponential contraction).
Outside the contractive regime, the picture remains incomplete. In Section \ref{gcbunique}, we explore further links between (non-)uniqueness and GCB.
In particular, we show that if a PCA has a translation-invariant Gibbs measure  $\mu $ that is both stationary and satisfies  GCB, then  $\mu $ must be the unique stationary
measure.

\section{Setting and basic notions}\label{pcasect}

\subsection{Configuration spaces and function spaces}

We consider two state spaces:
\begin{itemize}
\item[(i)]
 \( \Omega_S = \{-1,+1\}^{\mathbb{Z}^d} \), representing spatial (Ising) spin configurations,
\item[(ii)]
 \( \Omega_{ST} = \{-1,+1\}^{\mathbb{Z}^d \times \mathbb{N}} \), representing space-time spin configurations, where we adopt the convention \( \mathbb{N} = \{0,1,2,\ldots\} \) for discrete time.
\end{itemize}
For \( \sigma \in \Omega_S \), the spin value at lattice site \( x \in \mathbb{Z}^d \) is denoted by \( \sigma_x \in \{-1,1\} \). Similarly, for \( \sigma \in \Omega_{ST} \), the spin value at site \( x \in \mathbb{Z}^d \) at time \( n \) is denoted by \( \sigma_{x,n} \in \{-1,1\} \).

When stating definitions and results that apply to both \( \Omega_S \) and \( \Omega_{ST} \), we will simply use the notation \( \Omega \).

We write \( \Lambda \Subset \mathbb{Z}^d \) to mean that \(\Lambda\) is a finite subset of $\mathbb{Z}^d$.
For \( \Lambda \Subset \mathbb{Z}^d \), we denote by \( \mathcal{F}\!_\Lambda \) the \(\sigma\)-algebra generated by the spin evaluations \( \sigma \mapsto \sigma_x \) for \( x \in \Lambda \), with analogous notation for \( \Lambda \Subset \mathbb{Z}^d \times \mathbb{N} \).

Finally, for a function \( f: \Omega \to \mathbb{R} \) and \( \Lambda \Subset \mathbb{Z}^d \), we write \( f \in \mathcal{F}\!_\Lambda \) if \( f \) is measurable
with respect to \( \mathcal{F}\!_\Lambda \).

A probability measure \( \mu \) on \( \Omega_S \) is a probability measure defined on the Borel \(\sigma\)-algebra of \( \Omega_S \), where \( \Omega_S \) is endowed with the standard product topology, making it a compact metric space. The same definition applies to probability measures on \( \Omega_{ST} \).

We denote by \( \mu_\Lambda \) the restriction of \( \mu \) to the \(\sigma\)-algebra \( \mathcal{F}\!_\Lambda \). With a slight abuse of notation, we use \( \sigma_\Lambda \) to refer both to the restriction of a configuration \( \sigma \) to \( \Lambda \) and to the cylinder set in \( \mathcal{F}\!_\Lambda \) consisting of all configurations whose restriction to \( \Lambda \) matches that of \( \sigma \).
For \( \sigma \in \Omega_S \), we denote by \( \tau_x \sigma \) the shifted (or translated) configuration, defined by
\[
(\tau_x \sigma)_y = \sigma_{x+y}, \quad \text{for all } x,y \in \mathbb{Z}^d.
\]
A real-valued function \( f \) defined on either \( \Omega_S \) or \( \Omega_{ST} \) is called {\em local} if it depends only on a finite number of coordinates, that is, there exists
$\Lambda\Subset \Z^d$ such that \(f\in\mathcal{F}\!_\Lambda \). (When $\Lambda$ is the smallest such set, it is called the dependence set of $f$.)

By the Stone-Weierstrass theorem, the set of local functions is dense in the Banach space of continuous functions \( \mathcal{C}(\Omega) \) equipped with the supremum norm.

For a function \( f: \Omega \to \mathbb{R} \), we denote by \( \tau_x f \) the function defined by
\[
\tau_x f(\eta) = f(\tau_x \eta).
\]
The ``discrete gradient'' of a function $f:\Omega_S\to\R$ is defined as
\[
\nabla\!_x f (\si)= f(\si^{\sx})- f(\si)
\]
where $\si^{\sx}$ denotes the configuration obtained from $\si$ by flipping the spin at lattice site $x\in\Zd$.
We further define
\[
\delta\!_x f=\sup_{\si\in\Omega} \nabla\!_x f(\si)
\quad\text{and}\quad \delta f =(\delta_x f)_{x\in\Zd}.
\]
We can think of $\delta_x f$ as the oscillation of the function $f$ at site $x$.
For a local function we call its dependence set the finite set of those $x\in\Zd$ for which $\delta_x f\not= 0$.

For $p\geq  1$, we define the $\ell^p(\Zd)$ norm of $\delta f$ by
\[
\|\delta f\|_p^p= \sum_{x\in\Zd} (\delta_x f)^p.
\]
\subsection{Relative entropy}
For two probability measures $\mu,\nu$ and $\Lambda\Subset\Zd$ we define their relative entropy
in volume $\la\Subset\Zd$ via
\[
 \ent_{\Lambda} (\mu|\nu)=
\begin{cases}
\sum_{\si_\la\in\{-1,+1\}^\La} \mu(\si_\la) \log\frac{ \mu(\si_\la)}{ \nu(\si_\la)} & \text{if}\ \mu_\la\ll \nu_\la,
\\
\infty & \text{otherwise}.
\end{cases}
\]
For the sequence of cubes $C_n= [-n,n]^d\cap\Zd$, $n\geq 1$, we define the lower relative entropy density
\[
\liminf_{n\to\infty}\frac{\ent_{C_n} (\mu|\nu)}{|C_n|} = \ent_* (\mu|\nu).
\]
In case the limit exists, we write
\[
\lim_{n\to\infty}\frac{\ent_{C_n} (\mu|\nu)}{|C_n|} = s (\mu|\nu).
\]
for the relative entropy density.
This limit exists, for instance, when \( \nu \) is a Gibbs measure with respect to a translation-invariant, absolutely summable potential (see \cite{geo}, Chapter 15). More generally, it holds when $\nu$ is a translation-invariant probability measure satisfying the so-called asymptotic decoupling property (see \cite{pfister}).

\subsection{Gaussian concentration bound (GCB)}

We can now define what we mean by a Gaussian concentration for a given probability measure on $\Omega$.
Remember that $\Omega$ means either $\Omega_S$ or $\Omega_{S\!T}$.

\bd[Gaussian Concentration Bound]
\leavevmode\\
A probability measure $\mu$ on $\Omega$ is said to satisfy the Gaussian concentration bound with constant
$C>0$, abbreviated $\gcb{C}$, if for every local function $f:\Omega\to\R$, we have
\be\label{gcb}
\log\int \e^{f-\int f \dd\mu} \dd\mu\leq \frac{C}{2} \|\delta f\|_2^2.
\ee
\ed

For non-local continuous functions, inequality \eqref{gcb} is meaningful only when $\|\delta f\|_2<+\infty$. It is immediate from
\eqref{gcb} that the inequality can be extended to all uniform limits of local functions ({\em i.e.}, continuous functions) with $\|\delta f\|_2<+\infty$.
\br\label{prodrem}
\leavevmode
\ben
\item
If $\mu$ is a Dirac measure concentrating at some configuration $\si\in\Omega$, then \eqref{gcb} holds for all $C>0$. In fact in this case we can take $C=0$.
Conversely, since we always have $\int \e^{f-\int f \dd\mu} \dd\mu\geq 1$ (by Jensen inequality), if \eqref{gcb} holds for $C=0$, then necessarily $\mu=\delta_\si$ is a Dirac measure concentrating on a single configuration $\si$.
\item A product measure satisfies $\gcb{c}$ with a constant $c$ not depending on the marginals, see \cite[Theorem 6.2, page 161]{blm}.  One can choose $c=1/4$.
\een
\er
We have the following equivalent reformulation of \eqref{gcb} as a tail estimate, which follows from standard arguments starting from the exponential Chebychev inequality, see e.g. \cite[Proposition 2.5.2]{vershynin}.
\bp
A probability measure $\mu$ on $\Omega$ satisfies $\gcb{C}$ if and only if for all continuous $f:\Omega\to\R$ such that $\|\delta f\|_2<+\infty$, we have
\[
\mu\Big(f-\scaleto{\int}{18pt} f \dd\mu\geq u\Big)\leq \e^{-\frac{u^2}{2C \|\delta f\|_2^2}},\;\forall u>0.
\]
\ep
We recall a result from \cite{cr,moles} which we will use in Section \ref{bamba} below.
\bt\label{edithm}
If $\mu$ is a translation-invariant probability measure which satisfies $\gcb{C}$ then, for any translation-invariant  probability measure $\nu$ we have
\[
\ent_*(\nu|\mu)\geq \frac{\bar{d}(\nu,\mu)^2}{2C},
\]
where $\bar{d}(\cdot,\cdot)$ is the Dobrushin distance \cite{dobrushin}
defined on the set of translation-invariant probability measures by
\[
\bar{d}(\nu,\mu)= \sup\left\{\int f \dd\mu- \int f \dd\nu: \|\delta f\|_1\leq 1\right\}.
\]
In particular $\nu\not=\mu$ implies $\ent_*(\nu|\mu)>0$.
\et

\begin{remark}
The Dobrushin distance coincides with the so-called ``d-bar'' distance or Ornstein distance used in ergodic theory \cite{collet,shields}.
\end{remark}

\subsection{Probabilistic cellular automata (PCA)}\label{pcag}

We consider discrete-time  Markov processes on the configuration space $\Omega_S$.
We denote by $\si_n$ the configuration at time $n$, and $\si_{x,n}$ the value of the configuration at time $n$ at lattice site $x\in\Zd$.
Hence $(\si_n)_{n\in\N}\in \Omega_{S\!T}$.

We use the notation and framework of \cite{lms}. Conditional on $\eta$ at time $n\in\N$, the configuration one step later is distributed as a product measure $p(\dd\si|\eta)$ with marginals at site $x$ given by
\[
p_x(\eta):= p\big(\si_{x,n+1}= +1\, \big| \,\si_{n}=\eta\big)= \frac12\big(1 + \eta_x h_x (\si_{n})\big).
\]
Here $|h_x(\si)|\leq 1$ is given by
\be\label{hpara}
h_x (\eta)= \sum_{A\Subset \Zd} r\!_{\scaleto{A}{4.5pt}}\, \eta_{\scaleto{A}{4.5pt}+x},
\ee
where $\eta_{\scaleto{A}{4.5pt}}=\prod_{y\in A} \eta_y$. For the coefficients $r_A\in\R$ we assume moreover that
\be\label{bipo}
\sum_{A\Subset\Zd} |r\!_{\scaleto{A}{4.5pt}}|<\infty,
\ee
which immediately implies the absolute convergence of the sums appearing in \eqref{hpara}.
The $r\!_{\scaleto{A}{4.5pt}}$'s may be viewed as the parameters that govern the PCA.
We say that the PCA is of finite range $R$ if $r_A=0$ whenever the diameter of the set $A$ exceeds $R$.

We define the transition operator corresponding to a PCA as the operator acting on continuous functions defined by
\[
Pf(\eta)= \int f(\si)\, p(\dd\si|\eta)= \E\big(f(\si_1)|\si_0=\eta\big),
\]
where $f:\Omega_S\to\R$ is a continuous function. Notice that by \eqref{bipo}, $Pf$ is continuous whenever $f$ is continuous.
The corresponding operator acting on probability measures is then as usual defined via $\mu P (\dd \si) = \int p(\dd\si|\eta)\, \mu(\dd\eta)$,
and satisfies $\int Pf \dd\mu = \int f \dd\mu P$.
The form \eqref{hpara} of the functions $h_x$ ensure that the process is translation invariant, {\em i.e.}, the transition operator $P$ commutes with translations $\big\{\tau_x, x\in\Zd\big\}$.

We say that a probability measure \( \mu \) on \( \Omega_S \) is stationary if it satisfies \( \mu P = \mu \). By the compactness of the space of probability measures on \( \Omega_S \) (equipped with the topology of weak convergence), the set of stationary measures forms a non-empty, compact, and convex set, whose extreme points are known as ergodic measures.
In this work, we primarily focus on stationary and translation-invariant probability measures which also form a non-empty compact convex set. We denote by \( \mathcal{I}_\tau \) the set of such measures, and by \( \mathcal{I} \supset \mathcal{I}_\tau \) the  set of all invariant measures.

When $\nu\in \caI_\tau$ is a stationary translation invariant measure, the corresponding path space measure $\pee\!_\nu$ is a stationary process and can thus be extended to negative times, yielding a stationary and translation invariant space-time measure on
$\Omega_{d+1}= \{-1,1\}^{\Z^{d+1}}$. As proved in \cite{ms1}, this space-time measure is a translation invariant Gibbs measure for a Hamiltonian derived from the transition probabilities. An important consequence is that for $\nu_1,\nu_2\in\caI_\tau$, the corresponding (two-sided) path space measures $\pee\!_{\nu_1}, \pee\!_{\nu_2}$ have zero relative entropy density.
Another consequence is that the stationary measure $\nu$ is therefore the restriction of the space-time measure Gibbs $\pee_\nu$.
The space-time Gibbs measures of a PCA are not ``generic'' Gibbs measures due to the normalization condition of the transition probabilities \cite[Section 1.5]{lms}.


We denote by $\pee\!_\mu$ the path-space measure of the process $\{\eta_n,n\in\N\}$ starting from $\eta_0$ distributed according to $\mu$.
For a local function $h: \Omega_{S\!T}\to\R$, we define similarly, as in Section \ref{pcasect}, the oscillation $\delta_{x,n} h$, and the Gaussian concentration bound
where now $\|\delta h\|_2^2=\sum_{(x,n)\in\Zd\times\N} (\delta_{x,n} h)^2$.
\subsection{Contraction properties}
The following lemma is taken from \cite[Lemma 2, p. 136)]{lms}.
\bl\label{contlem}
Let $f:\Omega_S\to\R$ be a local function.
Then
\be\label{maeslem}
\delta_x (Pf) \leq \sum_{y\in\Zd} \sum_{A\ni x-y} |r\!_A|\, \delta_y (f),\;\forall x\in\Zd.
\ee
\el
The following corollary gives an estimate of the $\ell^2$-norm of $\delta(Pf)$.
\bc\label{contcol}
Define
$\psi(x)= \sum_{A\ni x} |r\!_A|$, $x\in\Zd$.
Let $f:\Omega_S\to\R$ be a local function. Then we have the estimate
\be\label{twonorm}
\|\delta (Pf)\|_2^2\leq \|\psi\|_1^2 \,\|\delta f\|_2^2,
\ee
where
\[
\|\psi\|_1= \sum_{x\in\Zd} \psi(x)= \sum_{A\Subset\Zd} |A||r\!_A|.
\]
As a consequence, we have
\[
\|\delta (P^n f)\|_2^2
\leq 
\|\psi\|_1^{2n} \|\delta f\|_2^2,\; \forall n\geq 1.
\]
\ec
\bpr
Notice that by the definition of $\psi$, the right-hand side of \eqref{maeslem} reads
$\sum_{y\in\Zd} \psi(x-y)\, \delta_y f= (\psi*\delta f)_x$. Because $\delta_x (Pf)\geq 0$ for all $x$,  we then have the estimate
\[
\|\delta (Pf)\|_2^2=\sum_{x\in\Zd}(\delta_x(Pf))^2\leq \sum_{x\in\Zd}(\psi*\delta f)_x^2= \|\psi*\delta f\|_2^2\leq \|\psi\|_1^2\, \|f\|_2^2\,,
\]
where in the last step we used Young's convolution inequality. We then obtain
\[
\|\psi\|_1=\sum_{x\in\Zd}\psi(x)= \sum_{x\in\Zd}\sum_{A\ni x} |r_A|
= \sum_{A\Subset \Zd} |A|\, |r_A|.
\]
The last statement follows by iteration:
\[
\|\delta (P^n f)\|_2^2
\leq 
\| \psi\|_1^2 \|\delta (P^{n-1} f)\|_2^2\leq\cdots\leq 
\|\psi\|_1^{2n} \|\delta f\|_2^2
\]
\epr
For later convenience we abbreviate
\be\label{kappa}
\kappa= \|\psi\|^2_1=\bigg(\,\scaleto{\sum_{A\Subset\Zd}}{22pt} |A|\, |r_A|\bigg)^2,
\ee
which by \eqref{twonorm} can be thought of as a contraction coefficient.
When $\kappa<1$ we say that the PCA is contractive, which implies uniqueness of the stationary measure $\mu$ and exponentially fast convergence to $\mu$ in the course of time.

\medskip



\section{GCB and time-evolution under a PCA}\label{pcaev}

\subsection{Propagation of GCB}
The first result shows that if we start from a probability measure satisfying $\gcb{C}$, then this property is conserved under the PCA evolution, with the associated constant evolving in time in a controlled manner, determined by the contraction coefficient $\kappa$ in \eqref{kappa}.

\bt[GCB is conserved under PCA evolution]\label{conservation}
\leavevmode\\
Let $\mu$ be a probability measure on $\Omega_S$ satisfying $\gcb{C}$, and $P$ the transition operator of a PCA, as
in Section \ref{pcasect}. Then we have the following.
\begin{enumerate}[\textup{(}a\textup{)}]
\item
$\mu P$ satisfies $\gcb{C_1}$ with
\[
C_1 = c+ C\kappa,
\]
where $c$ denotes the constant of the Gaussian concentration bound for a product measure (see Remark \ref{prodrem}, item 2). In particular we can choose $c=1/4$. 
\item
As a consequence, for all $n\in\N$,  $\mu P^n$ satisfies $\gcb{C_n}$ with
\[
C_n=
\begin{cases}
c\frac{\,\,1-\kappa^n}{1-\kappa} + C\kappa^n & \text{if}\;\; \kappa\not=1,
\\
c\,n+ C  & \text{if}\ \kappa=1.
\end{cases}
\]
\item
If $\kappa <1$, then $\mu P^n\to\nu$ where $\nu$, the unique stationary probability measure, satisfies $\gcb{C_\infty}$ with
\[
C_\infty= \frac{c}{1-\kappa}.
\]
\item Noninteracting case. When $h_x(\si_n)= h_x(\si_{n,x})$, we say that the PCA is non-interacting. In that case, GCB is conserved, and the unique stationary measure is a product measure, hence also satisfies GCB.
\end{enumerate}
\et
\bpr
Let $f:\Omega_S\to\R$ be a local function.
We have
\begin{align}\label{fipo}
\MoveEqLeft\int \mu P (\dd\si) \e^{f(\si)- \int f(\si)\mu P (\dd\sigma) }
\\
&=\int\mu(\dd\si) \int p(\dd\eta|\si) \e^{f(\eta) - \int f(\eta) p(\dd\eta|\si)} \e^{\int f(\eta) p(\dd\eta|\si)- \iint f(\eta) p(\dd\eta|\si) \mu(\dd\si)}.
\nonumber
\end{align}
Now we use that $p(d\eta|\si) $ is a product measure and as a consequence of $\gcb{c}$ we have, uniformly in $\si$:
\be\label{bala}
\int p(\dd\eta|\si) \e^{f(\eta) - \int f(\eta) \, p(\dd\eta|\si)}\leq \exp\left({\tfrac{c}{2}\|\delta f\|_2^2}\right).
\ee
Next, using the fact that $\mu$ satisfies $\gcb{C}$, combined with Corollary \ref{contcol} and \eqref{kappa}, shields
\begin{align}
\nonumber
\MoveEqLeft[4] \int \mu(\dd\si) \e^{\int f(\eta) p(\dd\eta|\si)- \iint f(\eta) p(\dd\eta|\si) \mu(\dd\si)}\\
&=\int \mu(\dd\si) \exp\left({ Pf(\sigma) - \int Pf(\si) \mu (\dd\si)}\right)
\nonumber\\
&\leq  \exp\left({\tfrac{C}{2} \|\delta (P f)\|_2^2}\right)
\nonumber\\
&\leq  \exp\left({\tfrac12 C \kappa\|\delta ( f)\|_2^2}\right).
\label{flipo}
\end{align}
Combining \eqref{bala} and \eqref{flipo}, we find
\[
\int \mu P (\dd\si) \e^{f(\si)- \int \mu P (\dd\sigma) f(\si)} \leq 
\exp\left({ \tfrac12(c+ C\kappa)\|\delta f\|_2^2}\right),
\]
which proves the first statement of the theorem.
The second statement follows by iteration, and the third one by taking the limit $n\to\infty$ which is possible when $\kappa<1$.
In case the PCA is non-interacting, a product measure is mapped to a product measure. Moreover, by the non-interacting property, we have for all $k\in \N$
and $x\in\Zd$ that $(\delta P^k f)_x\leq \delta_x f$ and $\delta_\si P^k$ is a product measure for all $k\in\N$.
Therefore, in that setting we can repeat \eqref{fipo}, \eqref{bala}, \eqref{flipo}, replacing $P$ by $P^k$ to conclude that $\mu P^k$ satisfies GCB with constant $C_k= C + c$.
\epr

\subsection{GCB of the space-time measure}

Our second main result shows that for a contractive PCA, {\em i.e.}, when $\kappa<1$, a space-time measure $\pee\!_\mu$ satisfies $\gcb{C}$ if
and only if the starting measure $\mu$ satisfies $\gcb{C'}$. One direction of this equivalence is obvious, because the Gaussian concentration bound
for $\pee\!_\mu$ immediately implies the Gaussian concentration bound for the single-time (at time zero) marginal $\mu$. The other implication is less obvious and shows that space-time Gibbs measures associated to PCA's are special ({\em i.e.}, not like generic Gibbs measures, see Section 4.1 below).

\bt\label{spacetime}
Consider a contractive PCA with transition operator $P$, that is, suppose that $\kappa<1$ (see \eqref{kappa} for the definition).
Let $\mu$ denote a probability measure on $\Omega_S$ and let $\pee\!_\mu$ denote the corresponding path-space measure.
Assume that $\mu$ satisfies $\gcb{C}$. Then there exists $C'$ such that  $\pee\!_\mu$ satisfies $\gcb{C'}$.
As a consequence, for the unique stationary measure $\nu$, the corresponding path-space measure
$\pee\!_\nu$ satisfies $\gcb{C'}$.
\et
\bpr
The consequence follows immediately from the fact that in the contractive regime $\kappa<1$, the
unique stationary measure satisfies $\gcb{C}$ for some $C>0$ by Theorem
\ref{conservation}.
We start by considering a local function $f(\si_0,\si_1)$ depending on the space-time configuration at times zero and one and estimate its centered moment generating function.
This example will make clear how to deal with a local function depending on the space-time configuration at times zero up to $n$.
We first fix some notation.
We define for a function $f: \Omega_{S\!T}\to\R$ its oscillation vector at time $n\in\N$ by
\[
\delta_{x}^n f= \sup_{\si\in \Omega_{S\!T}} \big(f(\si^{\scaleto{(x,n)}{6pt}})- f(\si)\big).
\]
We think of this as a vector indexed by $n$, and with $x$ the ``running'' index.
We have obviously
\[
\|\delta f\|_2^2 = \sum_{x\in\Zd,\, n\in\N} (\delta_{x}^{n} f)^2= \sum_{n\in\N} \|\delta^n f\|_2^2,
\]
where, with small abuse of notation,
\[
\|\delta^n f\|_2^2= \sum_{x\in\Zd} (\delta_{x}^{n} f)^2
\]
denotes the $\ell^2$-norm of the vector $\delta^n f$.

Let $f:\Omega_{S\!T}\to\R$ be of the form $ f(\si_0, \si_1)$. Then we write
\begin{align}\label{barastra}
\MoveEqLeft[4]\int \exp \left( f(\si_0, \si_1)- \int f \dd\pee_\mu\right) \dd\pee\!_\mu(\si_0, \si_1)=
\nonumber\\
& \hspace{-1.5cm} \iint \mu(\dd\si_0)\, p(\dd\si_1|\si_0)\left[ \exp \left(f(\si_0, \si_1)-  \int p(\dd\si_1|\si_0)\, f(\si_0, \si_1)\right)\right.
\nonumber\\
& \hspace{-1cm}\left.\exp \left( \int p(\dd\si_1|\si_0)\, f(\si_0, \si_1)-\iint\mu(\dd\si_0) \, p(\dd\si_1|\si_0)\, f(\si_0, \si_1)\right)\right].
\end{align}
Using that $p(\dd\si_1|\si_0)$ is a product measure which therefore satisfies $\gcb{c}$, we start then by estimating
\begin{align}\label{inner}
\MoveEqLeft[4] \int p(\dd\si_1|\si_0)\,\exp \left(f(\si_0, \si_1)-  \int p(\dd\si_1|\si_0)\, f(\si_0, \si_1)\right)
\nonumber\\
&\leq \exp \big(\tfrac12 c\|\delta f_{\si_0}\|_2^2\big)\,,
\end{align}
where $f_{\si_0}:\Omega_S\to\R$ is the function defined by $f_{\si_0} (\si)= f(\si_0, \si)$.
As a consequence, we have for all $\si_0\in\Omega_S$ the estimate $\delta_x f_{\si_0}\leq \delta_{x}^1 f$, and therefore,
\[
\|\delta f_{\si_0}\|_2^2\leq \|\delta^1 f\|_2^2\,.
\]
Next, using the fact that $\mu$ satisfies $\gcb{C}$ we tackle the second factor in \eqref{barastra} as follows:
\begin{align}\label{outer}
\MoveEqLeft \int \mu(\dd\si_0)\exp \left( \int p(\dd\si_1|\si_0) f(\si_0, \si_1)-\iint\mu(\dd\si_0)\, p(\dd\si_1|\si_0) f(\si_0, \si_1)\right)
\nonumber\\
&=
\int \mu(\dd\si_0) \exp \left(F(\si_0)-\int F(\si_0)\mu(\dd\si_0)\right)
\nonumber\\
&\leq
\exp \left(\tfrac{C}{2} \|\delta F\|_2^2\right),
\end{align}
where we denoted  $F(\si)= \int p(\dd\si_1|\si) f(\si,\si_1)$. Then we are left with estimating $\|\delta F\|_2^2$. We have
\begin{align*}
&F(\si^{\sx})- F(\si)\\
&= \int p(\dd\si_1|\si^{\sx}) \left(f(\si^{\sx},\si_1)- f(\si, \si_1)\right)\\
& \qquad + \int p(\dd\si_1|\si^{\sx}) f(\si,\si_1)- \int p(\dd\si_1|\si) f(\si,\si_1)\\
&\leq \delta_{x,0} f + \big(P (f_\si))[\si^{\sx}]- (P(f_\si)\big) [\si],
\end{align*}
where $f_\si:\Omega\to\R: \eta\mapsto f(\si,\eta)$, and as a consequence,
\[
Pf_\si (\eta)= \int f_\si(\xi) \,p(\dd\xi|\eta)= \int f(\si, \xi)\, p(\dd\xi|\eta).
\]
Therefore, using Lemma \ref{contlem}, we obtain
\be\label{peest}
(P (f_\si))\big(\si^{\sx}\big)- (P(f_\si)) (\si)\leq (\psi*\delta f_\si)_x\leq (\psi* \delta^1 f)_x.
\ee
Combining then \eqref{inner}, \eqref{outer}, \eqref{peest}, we obtain
\[
\log \int \exp \left( f(\si_0, \si_1)- \int f \dd\pee\!_\mu\right) \dd\pee\!_\mu \leq \tfrac{c}{2} \, \big\|\delta^1 f \big\|_2^2+  
\tfrac{C}{2}\Big(\big\|\delta^0 f + \psi* \delta^1 f \big\|_2^2\Big).
\]
Considering now a function $f(\si_0,\si_1,\ldots,\si_n)$, that is, a function of the configuration at times $0$ up to $ n$, we can iteratively proceed as before by writing
\begin{align*}
\MoveEqLeft \int \exp \left( f(\si_0,\ldots,\si_n)- \int f \dd\pee\!_\mu\right) \dd\pee\!_\mu\\
&=\int \mu(\dd\si_0) \int p(\dd\si_1|\si_0)\cdots \int p(\dd\si_n|\si_{n-1})\\
&\qquad \Bigg[\exp\Big(f- \int p(\dd\si_n|\si_{n-1})f(\si_0,\ldots, \si_{n-1},\si_n)\Big)\\
&\;\;\qquad\exp\Bigg(\int p(\dd\si_n|\si_{n-1})f(\si_0,\ldots, \si_{n-1},\si_n)\\
& \quad\qquad\qquad -\int p(\dd\si_{n-1}|\si_{n-2})\, p(\dd\si_n|\si_{n-1})f(\si_0,\ldots, \si_{n-1},\si_n)\Bigg)\\
&\qquad\quad \cdots \exp\left(\int f_{\si_0} \dd\mu-\int f\dd\pee\!_\mu\right)\Bigg].
\end{align*}
Then we estimate the different factors, from inner to outer, by using $\gcb{c}$ for product measures, and for the final (outer) integral
$\gcb{C}$ for the measure $\mu$.
This gives the following estimate
\begin{align}\label{baster}
\MoveEqLeft 2\log\int \exp \left( f(\si_0, \si_1,\ldots,\si_n)- \int f \dd\pee\!_\mu\right) \dd\pee_\mu(\si_0, \si_1,\ldots,\si_n)
\nonumber\\
&\leq
c\, \|\delta^n f\|_2^2
\nonumber\\
&\quad +
c\, \big\|\delta^{n-1}f + \psi* \delta^{n} f\big\|_2^2
\nonumber\\
& \quad\;\, \vdots \nonumber \\
&\quad +
c\, \big\| \delta^1 f + \psi* \delta^2 f + \cdots + (\psi^*)^{n-1}*\delta^n f \big\|_2^2
\nonumber\\
&\quad +  C\,\big\| \delta^0 f + \psi* \delta^1 f + \cdots + (\psi^*)^{n}*\delta^n f \big\|_2^2,
\end{align}
where $(\psi^*)^k$ denotes the $k$-fold convolution power of $\psi$.
We are now left with estimating the rhs of \eqref{baster}.

Now recall the notation \eqref{kappa}, then  we have, by iterative application of Young's inequality,
\begin{equation}\label{kopu}
\|(\psi*)^k*(\delta^j f)\|_2\leq \|\psi\|_1^k\, \|(\delta^j f)\|_2 =\kappa^{k/2}\, \|(\delta^j f)\|_2
\end{equation}
As as consequence, combining \eqref{kopu} with by Minkowski inequality
$\|f+g\|_2\leq \|f\|_2+\|g\|_2$, we can estimate the terms appearing in the rhs of \eqref{baster} as follows
\begin{align*}
& \big\|\delta^{n-1}f + \psi* \delta^{n} f\big\|_2^2 \; \leq \; \big(\big\|\delta^{n-1} f\|_2 + \sqrt{\kappa}\, \|\delta^n f\big\|_2\big)^2\\
& \qquad\qquad \vdots \\
& \big\| \delta^1 f + \psi* \delta^2 f + \cdots + (\psi^*)^{n-1}*\delta^n f \big\|_2^2\\
& \qquad\qquad\qquad\qquad\;\;\;\leq \Big( \big\| \delta^1 f \big\|_2 + \sqrt{\kappa}\,\big\| \delta^2 f \big\|_2 + \cdots + (\sqrt{\kappa})^{n-1} \big\|\delta^n f \big\|_2\Big)^2.
\end{align*}
As a consequence, the sum in the rhs of \eqref{baster} can be estimated by
$\| A v\|_{\R^{n+1}}^2$ where $\|\cdot \|_{\R^{n+1}}$ denotes Euclidean norm in $\R^{n+1}$ and where
$v$ is the vector with elements $v_i= \|\delta^{i}f\|_2$, for $i=0,1,\ldots, n$, and $A$ is the lower triangular $(n+1)\times (n+1)$ matrix whose elements are given by
\begin{align}\label{mat}
A_{i,n-j} &= \sqrt{c} \,\sqrt{\kappa}\,^{n-j}\ \  0\leq i\leq n-1,\ \ 0\leq j\leq i
\nonumber\\
A_{n, n-j} &= \sqrt{C} \,\sqrt{\kappa}\,^{n-j}\ \ \   0\leq k\leq n
\nonumber\\
A_{k,l} &= 0 \ \text{otherwise}.
\end{align}
The Euclidean norm $\|Av\|_2^2$ can therefore be estimated by $\|A\|_2^2\, \|v\|^2$. The matrix norm of $A$ can in turn be estimated by
\[
\|A\|_2^2 \leq \|A\|_\infty \|A\|_1,
\]
where
\[
\|A\|_\infty= \sup_i \sum_j |A_{i,j}|, \ \|A\|_1= \sup_j \sum_i |A_{i,j}|.
\]
Using  \eqref{mat} and the fact that $\kappa<1$, which implies that, for all $j$,
$\sum_{i=0}^j (\sqrt{\kappa})^{i}\leq \frac{1}{1-\sqrt{\kappa}}$, we obtain
\[
\|A\|_2\leq \frac{\sqrt{c}}{1-\sqrt{\kappa}}\vee \frac{\sqrt{C}}{1-\sqrt{\kappa}},
\]
whereas
\[
\| v\|^2= \sum_{k=0}^{n} \|\delta^k f\|_2^2 =  \|\delta f\|_2^2.
\]
As a consequence we obtain as a final estimate for the sum in the rhs of \eqref{baster}
\[
\|A\|_2^2\, \|v\|^2
\leq \left(\frac{\sqrt{c}}{1-\sqrt{\kappa}}\vee \frac{\sqrt{C}}{1-\sqrt{\kappa}}\right)^2\,\|\delta f\|_2^2\,,
\]
which shows that we have the Gaussian concentration bound with constant
\[
C'=\left( \frac{\sqrt{c} \vee \sqrt{C}}{1-\sqrt{\kappa}}\right)^2.
\]
The proof of the theorem is complete.
\epr

\begin{remark}
In the proof of Theorem \ref{spacetime} we used the contractive property $\kappa<1$ only in the estimation of
the norm of the matrix $A$.
As a consequence, we have that, if $\mu$ satisfies $\gcb{C}$, then the joint distribution
of $\si_0,\ldots, \si_n$ (viewed as a measure on $\{-1,1\}^{\Zd\times \{0,\ldots,n\}}$) satisfies $\gcb{\widetilde{C}_n}$
with ${\widetilde{C}_n}= \|A\|^2_2$, with $A$ the matrix given in \eqref{mat}.
If $\kappa<1$, then the matrix norm can be bounded uniformly in $n$, and therefore in that case we indeed can
consider the ``limit'' $n\to\infty$, \em{i.e.}, conclude that $\pee\!_\mu$ satisfies $\gcb{C'}$, with $C'=\sup_n  \widetilde{C}_n$.
\end{remark}
\begin{remark}
Because we have proved that the equivalence ``$\mu$ satisfies GCB with $\pee\!_\mu$ satisfies GCB'' only in the contractive regime, we cannot exclude that there are examples where $\mu$ (a stationary measure) satisfies GCB, but $\pee\!_\mu$ does not satisfy GCB, and possibly there are other stationary measures. As we will see below, in that case we can conclude that at least one of these stationary measures is non-Gibbs.
\end{remark}

\subsection{Relaxation to equilibrium via GCB and contraction estimates}\label{bamba}

In the contractive regime $\kappa<1$, there is a unique stationary measure $\nu$, and
from any initial measures $\mu$ we have that $\mu P^n\to \nu $ exponentially fast.
In this section we show how to quantify the convergence to $\mu$, in terms of a distance
(giving a stronger convergence than weak convergence), under abstract locality and a 
contraction conditions.
We start by defining these locality and contraction conditions.
\bd
\leavevmode
\ben
\item
We then say that the PCA is local with propagation speed $a\in \N$ if for all $n,k\in\N$ and
$f\in \caF_{C_n}$, we have $P^k f\in \caF_{C_{n+ ak}}$. Here as before, $C_n= [-n,n]^d\cap\Zd$.
\item We say that the PCA is contractive with contraction speed $\psi_n:\Zd\to\R^+$ if for all $n\in\N$ we have
\be\label{contspeed}
(\delta (P^n f))\leq \psi_n*\delta f,
\ee
where the inequality has to be interpreted point-wise, \em{i.e.}, $(\delta P^n f)_x\leq (\psi_n*\delta f)_x$, for all $x\in\Zd$.
\een
\ed
The condition on locality is satisfied for finite range PCA's.

In order to formulate our result, we need the following pseudo-distances.
\bd
For $\la\Subset\Zd$ we define
\[
D_{\infty,\la}(\mu,\nu) = \sup \left\{ \int f \dd\mu-\int f \dd\nu: {f\in\caF_\la},  \|\delta f\|_1\leq 1\right\},
\]
\[
D_{2,\la} (\mu,\nu)= \sup\left\{\int f \dd\mu-\int f \dd\nu: {f\in\caF_\la},  \|\delta f\|_2\leq 1\right\}.
\]
\ed
Notice that $D_{\infty,\la}(\mu,\nu), D_{2,\la} (\mu,\nu)>0$ if and only if $\mu_\la \not=\nu_\la$.

The Cauchy-Schwarz inequality  immediately gives
$\|\delta f\|_1\leq \sqrt{|\la|} \|\delta f\|_2$ and, as a consequence,
\[
\frac{D_{2,\la}}{\sqrt{|\la|}}\leq D_{\infty,\la}.
\]

For $\nu,\mu$ translation-invariant probability measures we have (see \cite{collet})
\[
\lim_{\la\uparrow\Zd}D_{\infty,\la}(\nu,\mu)=\bar{d}(\nu,\mu).
\]
The $\bar{d}$-distance is well-known and widely used in ergodic theory, see e.g. \cite{shields}. Convergence in this distance is stronger than weak convergence and  conserves, for instance,  ergodicity, mixing, and the relative entropy density w.r.t.\ a translation-invariant Gibbs measure is a continuous function w.r.t.\ this metric.
Additionally, we prove in \cite{collet} that
\[
\lim_{n\to\infty}\frac{D_{2, C_n}}{\sqrt{|C_n|}} = \bar{d}(\nu,\mu).
\]
Then we have the following result.
\bt
Assume the PCA has contraction speed $\psi$, is local with propagation speed $a\in\N$, and has a stationary measure $\mu\in \caI_\tau$ which
satisfies $\gcb{C}$.
Assume additionally that $\mu$ has ``finite energy'', {\em i.e.}, there exists $\rho>0$ such that
\[
\mu(\si_\la)\geq \e^{-\rho|\la|},\;\forall \sigma\in\Omega_S, \la\Subset\Zd.
\]
Then, for every probability measure $\nu$, and for all $n,k\in\N$, we have 
\be\label{dinf}
D^2_{\infty,C_n} (\nu P^k,\mu )\leq  2C \ent_{C_{n+ak}} (\nu|\mu) \,\|\psi_k\|_2^2,
\ee
and
\be\label{2dist}
D^2_{2,C_n} (\nu P^k,\mu )\leq 2C {\ent_{C_{n+ak}} (\nu|\mu)} \, \|\psi_k\|_1^2.
\ee
As a consequence,  for all $n,k\in\N$,
\be\label{dinfconv}
D^2_{\infty,C_n} (\nu P^k,\mu )\leq 2C\rho\, {|C_{n+ak}|}\|\psi_k\|_2^2.
\ee
If $\nu$ is a translation-invariant probability measure, then for all $k\in\N$
\be\label{dbarconv}
\bar{d}\,^2(\nu P^k,\mu)\leq 2C\rho\, \|\psi_k\|_1^2.
\ee
As a consequence, if $\|\psi_k\|_1\to 0$ as $k\to \infty$, then $\nu P^k\to\mu$ in $\bar{d}$-distance.
\et
\bpr
First observe that by the finite energy condition, for every probability measure $\nu$ we have the upper bound
\be\label{finergb}
\frac{\ent_\la(\nu|\mu)}{|\la|}\leq \rho.
\ee
We then start with the variational formula for the relative entropy combined with the locality condition, the stationarity of $\mu$, and \eqref{contspeed} to obtain
\begin{align*}
\ent_{C_{n+ak}} (\nu|\mu) &\geq  \sup_{f\in \,\caF\!_{C_n}}\left( \int P^k f \dd\nu-\int P^k f \dd\mu - \log\int \e^{P^k f-\int f \dd\mu} \dd\mu\right)\\
&\geq \sup_{f\in \,\caF\!_{C_n}}\left( \int f \dd\nu P^k-\int  f \dd\mu - \frac{C}{2}\|\psi_k*\delta f\|_2^2\right).
\end{align*}
Now we use
\[
\|\psi_k*\delta f\|_2^2\leq \|\psi_k\|_2^2\, \|\delta f \|_1^2.
\]
Replacing then $f$ by $\lambda f$ and optimizing over $\lambda$ we obtain
\eqref{dinf}.
If on the other hand we use
\[
\|\psi_k*\delta f\|_2^2\leq \|\psi_k\|_1^2 \,\|\delta f \|_2^2,
\]
and proceed with the same optimization, then we arrive at \eqref{2dist}. The result
\eqref{dinfconv} then follows immediately via \eqref{finergb}, whereas
\eqref{dbarconv} follows by first dividing both sides of \eqref{2dist} by $|C_n|$ and then taking the limit $n\to\infty$, using once more \eqref{finergb}, combined with the fact that for all $k\in \N$ we have
$\tfrac{|C_{n+ak}|}{|C_n|}\to 1$ as $n\to\infty$.
\epr
\br
Notice that if $\kappa<1$ (cf. \eqref{kappa}) then $\|\psi_1\|<\kappa$ and by iteration of Young's inequality we have
\[
\|\psi_k*\delta f\|_2^2\leq \kappa^k \|\delta f\|_2^2\,.
\]
This then leads to exponential convergence in $\bar{d}$-distance, {\em i.e.},
\[
\bar{d}(\nu P^k,\mu)\leq 2C\rho\, \kappa^k.
\]
The exponential convergence in $\bar{d}$-distance for a class of PCA's and interacting particle systems was obtained in \cite{steiff} using coupling methods.
\er


\section{GCB and uniqueness}\label{gcbunique}

In \cite{moles} we provided a link between GCB and the uniqueness of transla\-tion-invariant Gibbs measures in the context of lattice spin systems with
translation-invariant potentials, see also \cite{cr} for a generalized setting beyond Ising spins.
Here we prove some connections between GCB and uniqueness in our current setting of PCA's for the space-time stationary measures and the stationary measures. Since we
use relative entropy density, in this section our results are necessarily restricted to translation invariant stationary measures.

\bt\label{uniprop}
If $\nu\in\caI_\tau$ is such that the corresponding path space measure $\pee\!_\nu$ satisfies $\gcb{C}$, then $\caI_\tau= \{\nu\}$. Moreover the measure $\nu$ is (temporally) mixing for the PCA, \em{i.e.}, for all $\nu$-centered local functions $f$ we have
\be\label{bongo}
\lim_{n\to\infty}\int f (P^n f) \dd\nu= 0.
\ee
\et
\bpr
If $\nu'\in\caI_\tau$ is another stationary translation-invariant measure, then the corresponding path space measure $\pee\!_{\nu'}$ has zero relative entropy density w.r.t. $\pee\!_\nu$, because both are Gibbs for the same translation-invariant potential (cf. \cite{gold}) . This fact, combined with the fact that $\pee\!_\nu$ satisfies $\gcb{C}$ by assumption, and using Theorem
\ref{edithm}, gives  $\pee\!_{\nu'}=\pee\!_{\nu}$.
As a consequence, by restriction to the time zero layer, we obtain $\nu=\nu'$.
The mixing property \eqref{bongo} follows from the space-time mixing property of $\pee\!_\nu$ which in turn follows from the GCB property of $\pee\!_\nu$ (see \cite[Proposition 2.2]{moles}).
\epr
\br
Within the framework of Theorem \ref{uniprop}, even if the stationary measure \( \nu \) is unique, we cannot directly conclude the convergence \( \mu P^n \to \nu \). Moreover, it remains unclear whether there exist examples of \( \pee\!_\nu \) satisfying \( \gcb{C} \) outside the contractive regime.
\er

The next result proves that certain forms of non-uniqueness are not compatible with GCB. In particular, starting from a ``low-temperature phase'' one cannot arrive at measure satisfying GCB in finite time.

We call the transition operator of a PCA non-degenerate if for all $\mu,\nu$ probability measures
on $\Omega_S$ we have $\mu\not=\nu$ implies $\mu P^n\not= \nu P^n$ for all $n\in\N$. This is satisfied, {\em e.g.},  when the range of  $P^n$ is uniformly dense in the set of continuous functions.

\bt\label{boprop}
\leavevmode
\ben
\item
Assume $P$ is non-degenerate. Let $\mu\not=\mu'$ denote two translation-invariant probability measures with relative entropy density $\ent(\mu|\mu')=0$.
Then for all $n\in\N$, $\mu' P^n$ does not satisfy $\gcb{C}$ for any $C\geq 0$.
\item Let $\nu,\nu'\in\caI_\tau$.  If $\nu$ satisfies $\gcb{C}$ then $\ent_*(\nu'|\nu)>0$. In particular, $\nu$ and $\nu'$ cannot be Gibbs measures for the same
translation-invariant potential.
\item Let $\mu$ be a translation-invariant Gibbs measure w.r.t. a translation invariant potential which is also stationary for the PCA. Assume $\mu$ satisfies $\gcb{C}$. Let $\nu$ be a translation-invariant probability measure with zero relative entropy production w.r.t. $\mu$, \em{i.e.}, such that
\be\label{zeroprod}
\ent(\nu|\mu)= \ent(\nu P|\mu),
\ee
then $\nu=\mu$. As a consequence $\caI_\tau=\{\mu\}$.
\een
\et
\bpr
For item $1$: $\ent(\mu|\mu')=0$ implies $\ent(\mu P^n|\mu' P^n)=0$ for all $n\in\N$, and therefore, if $\mu' P^n$ satisfies $\gcb{C}$, $\mu P^n= \mu' P^n$ which is impossible by the assumed non-degeneracy.

Item 2) From Theorem \ref{edithm} it follows that for $\nu'\not=\nu$, $\ent_*(\nu'|\nu)>0$. As a consequence $\nu$ and $\nu'$ cannot be Gibbs measures for the same translation-invariant potential, because this implies $\ent(\nu'|\nu)=0$.

For item 3: in \cite{daipra} it is proved that \eqref{zeroprod} implies that $\nu$ is a Gibbs measures for the same translation-invariant potential, and therefore, $\ent(\nu|\mu)=0$ which implies $\nu=\mu$ because $\mu$ satisfies $\gcb{C}$ by assumption. In particular when $\nu\in\caI_\tau$ we have \eqref{zeroprod}, and therefore $\nu=\mu$.
\epr

\subsection{ Two examples}
Finally, we provide  two examples.
The first example illustrates that outside the context of PCA's, it is possible to have two different ``space-time Gibbs measures'' and at the same time restrictions which satisfy GCB. The second example shows that in the context of PCA's one can have one invariant measure satisfying GCB, and another one which is both non-Gibbs and does not satisfy GCB.
\subsubsection*{Example 1: Restrictions in the low-temperature Ising model}

Item 3 of Theorem \ref{boprop} implies that if there exists a stationary measure \( \mu \) satisfying \( \gcb{C} \) and it is not unique -- more precisely, if there exists another stationary, translation-invariant measure \( \nu \) -- then \( \mu \) cannot be a Gibbs measure for a translation-invariant potential.

As discussed in Section \ref{pcag}, a stationary (translation-invariant) measure of a PCA is the restriction of a translation-invariant space-time Gibbs measure to the time-zero layer
\( \big\{(x,0): x \in \mathbb{Z}^d\big\} \) (or any other time layer). Consequently, if a stationary measure is both a Gibbs measure and satisfies \( \gcb{C} \), then it must be unique. The ability to conclude uniqueness from the Gibbs and GCB properties of a restriction is specific to PCA space-time Gibbs measures, as we will illustrate with the following example.

Consider the decimation of the plus and minus phases of the low-tempera\-ture Ising model in dimension \( d=2 \) and decimation parameter $b\geq 3$, i.e.,
the restrictions of the plus and minus phases to the sublattices
$b\Z^d$ with $b\geq 3$. Let \( \pee^+ \) and \( \pee^- \) denote these two translation-invariant measures. Their restrictions to the horizontal axis  denoted \( \nu^+ \) and \( \nu^- \), are Gibbs measures with exponentially decaying interactions and therefore satisfy \( \gcb{C} \) (see \cite{mrv}, \cite{mvdv}, and \cite{cckr} for the resulting GCB property). However, the full measures \( \pee^+ \) and \( \pee^- \) do not satisfy GCB, as their relative entropy satisfies \( \ent(\pee^+\vert\,\pee^-) = 0 \) and they are distinct measures, \( \pee^+ \neq \pee^- \). By Theorem \ref{boprop}, Item 3, such a situation cannot arise for space-time Gibbs measures of a PCA. In fact, in this example, both \( \pee^+ \) and \( \pee^- \) are non-Gibbsian measures \cite{enter}.

Finally, since \( \nu^+ \) satisfies \( \gcb{C} \), we also obtain \( \ent(\nu^+|\nu^-) > 0 \), meaning that \( \nu^+ \) and \( \nu^- \) cannot be Gibbs measures for the same potential.

\subsubsection*{Example 2: Cellular automata with totally asymmetric noise}

An important class of PCA's including well-known models such as Toom's model (in $d=2$), Stavskaya's model (in $d=1$) are studied in \cite{fertoom}. For the class of models considered in this work, the authors prove that in the ``high-noise'' regime, the unique invariant measure is $\delta_{+}$ (the Dirac measure concentrating on the all-plus configuration), and in the ``low-noise'' regime, there exist multiple invariant measures (including $\delta_{\pm}$). For each of these invariant measures $\mu$ one can prove that
\be\label{ferto}
\mu( +_{C_n})\geq \e^{-c n^{d-1}},
\ee
where $+_{C_n}$ denotes the event to see the all-plus configuration in the cube $C_n$. 
In \cite{poncelet} a complementary upper bound is proved.
From \eqref{ferto} we conclude that $s_*(\delta_+|\mu)=0$

If $\mu\in \caI_\tau$, $\mu\not=\delta_+$ satisfied $\gcb{C}$ then, by Theorem \ref{edithm} we would conclude
\[
\ent_*(\delta_+|\mu)>0,
\]
which contradicts \eqref{ferto}.
Therefore, in the regime of non-uniqueness for all these examples, from the set of invariant measures only the Dirac measure $\delta_+$ satisfies GCB, and all other invariant measures do not satisfy GCB.
Moreover, when starting from the all-minus configuration, in all these examples one converges
exponentially fast to an invariant measure $\mu^-$ which does not satisfy GCB. I.e., the path space measure starting from the measure $\delta_{-}$ (which trivially satisfies GCB) does not satisfy GCB, even if at every positive time $k>0$, the measure $\delta_ P^k$ does satisfy GCB (by Theorem \ref{conservation}).

\bigskip

\noindent \textbf{Acknowledgement.}
The authors thank Christian Maes for useful advise on PCA's and associated Gibbs measures.
The authors also thank Pierre Collet for useful discussions.



\begin{thebibliography}{99}

\bibitem{blm}
S. Boucheron, G. Lugosi, P. Massart.
{\em Concentration inequalities. A nonasymptotic theory of independence.} Oxford University Press, 2013.


\bibitem{ccr}
J.-R. Chazottes, P. Collet, F. Redig,  On concentration inequalities and their applications for Gibbs measures in lattice systems. J. Stat. Phys.  169  (2017),  no. 3, 504-546.



\bibitem{cckr}
J.-R. Chazottes, P. Collet,  C. K\"ulske, F.  Redig,  Concentration inequalities for random fields via coupling. Probab. Th. Rel. Fields  137  (2007),  no. 1-2, 201-225.

\bibitem{ccr2}
J.-R. Chazottes, P. Collet, F.  Redig, Evolution of concentration under lattice spin-flip dynamics. J. Statist. Phys. 184 (2021): 1-21.

\bibitem{moles}
J.-R. Chazottes, J. Moles, F. Redig, E. Ugalde,
Gaussian concentration and uniqueness of equilibrium states in lattice systems. J. Statist. Phys. 181 (2020), 2131-2149.


\bibitem{collet}
J.-R. Chazottes, P. Collet, F.  Redig.
Gaussian concentration, integral probability metrics and nonlinear Kantorovich-Rubinstein metrics for lattice systems, preprint 2025.

\bibitem{cr}
J.-R. Chazottes, F. Redig.
Relative entropy, Gaussian concentration and uniqueness of equilibrium states. Entropy, 24 (11) (2022), 1513.

\bibitem{daipra}
P. Dai Pra, P.-Y. Louis, S. Roelly. Stationary measures and phase transition for a class of probabilistic cellular automata. ESAIM: Probability and Statistics, 6, (2002), 89-104.

\bibitem{dobrushin}
R.L. Dobrushin. Prescribing a system of random variables by conditional distributions, Theory of Probability and Its Applications, 15 (1970), 458-486.

\bibitem{dobshlos}
R. L. Dobrushin, S.  Shlosman.
Completely analytical interactions: constructive description. J. Statist. Phys. 46 (1987): 983-1014.

\bibitem{dp}
D. P. Dubhashi, A. Panconesi.
{\em Concentration of Measure for the Analysis of Randomized Algorithms.}
Cambridge University Press, 2009.


\bibitem{enter}
A. C. D. van Enter, R. Fern\'{a}ndez, A. D. Sokal. Regularity properties and pathologies of position-space renormalization-group transformations: scope and limitations of Gibbsian theory. J. Statist. Phys. 72 (1993), no. 5-6, 879-1167.

\bibitem{fertoom}
R. Fern\'{a}ndez,  A. Toom.
Non-Gibbsianness of the invariant measures of non-reversible cellular automata with totally asymmetric noise. Ast\'{e}risque 287 (2003), 71-87.

\bibitem{geo}
H.-O. Georgii, {\em Gibbs measures and phase transitions}, second edition, De Gruyter, 2011.

\bibitem{kuelske}
C. K\"{u}lske,
Concentration inequalities for functions of Gibbs fields with application to diffraction and random Gibbs measures. Comm. Math. Phys.  239  (2003),  no. 1-2, 29-51.


\bibitem{gold}
S. Goldstein, R.  Kuik, J. L.  Lebowitz, and C. Maes.
From PCA's to equilibrium systems and back. Commun. Math. Phys. 125 (1989), 71-79.

\bibitem{lms}
J. Lebowitz,  C. Maes, and E.R. Speer, Statistical mechanics of probabilistic cellular automata. J. Stat. Phys. 59 (1990): 117-170.

\bibitem{ledoux}
M. Ledoux. {\em The concentration of measure phenomenon}. American Mathematical Soc., 2001.



\bibitem{mrv}
C. Maes, F. Redig, and A. Van Moffaert.  The restriction of the Ising model to a layer. J. Stat. Phys. 96  (1999), 69-107.

\bibitem{ms1}
C. Maes and S. Shlosman.  Ergodicity of Probabilistic Cellular Automata: A Constructive Criterion, Commun. Math. Phys. 135 (1991), 233-251.


\bibitem{mm}
J. Mairesse, I. Marcovici. Around probabilistic cellular automata, Theor. Comp. Sc. 559 (2014), 42-72.

\bibitem{mvdv}
J. L\"{o}rinczi, C. Maes, and K. V. Velde. Transformations of Gibbs measures. Probab. Th. and Rel. Fields 112 (1998), 121-147.

\bibitem{pfister}
C.-E. Pfister. Thermodynamical aspects of classical lattice systems. In and out of equilibrium (Mambucaba, 2000), 393-472, Progr. Probab., 51, Birkha\"user, 2002.

\bibitem{poncelet}
L. Ponselet. Phase transitions in probabilistic cellular automata. arXiv preprint arXiv:1312.3612 (2013).


\bibitem{shields}
P.  C. Shields. The ergodic theory of discrete sample paths. Vol. 13. American Mathematical Soc., 1996.

\bibitem{steiff}
J. E. Steif. Convergence to equilibrium and space-time bernoullicity for spin systems in the $M< \epsilon$ case.
Ergodic Theory and Dynamical Systems 11.3 (1991), 547-575.

\bibitem{stroock}
D. W. Stroock and  B. Zegarlinski. The equivalence of the logarithmic Sobolev inequality and the Dobrushin-Shlosman mixing condition.
Commun. Math. Phys. 144 (1992), 303-323.

\bibitem{vershynin}
R. Vershynin. {\em High-dimensional probability: An introduction with applications in data science}. Cambridge university press, 2018.

\bibitem{Wbook}
M. Wainwright.
{\em High-Dimensional Statistics: A Non-Asymptotic Viewpoint}.
(Cambridge Series in Statistical and Probabilistic Mathematics). Cambridge University Press, 2019.


\end{thebibliography}
\end{document}